\documentclass{article}
\usepackage{a4wide,amsmath,amsfonts,dsfont,latexsym,amssymb,amscd, graphicx,verbatim}
\usepackage{amsthm}

\newtheorem{theorem}{Theorem}[section]
\newtheorem{corollary}{Corollary}[section]
\newtheorem{lemma}{Lemma}[section]

\newcommand{\R}{{\mathbb R}}

\newcommand{\eps}{{\varepsilon }}
\newcommand{\tqed}{\text{$\blacksquare$}}

\newtheorem{thm}{Theorem}[section]

\newtheorem{prop}[thm]{Proposition}
\numberwithin{equation}{section}

\begin{document}
\title{Asymptotic expansion of the minimum covariance determinant estimators}
\author{Eric A.~Cator and Hendrik P.~Lopuha\"a\\
\\
\emph{Delft University of Technology}}
\date{\today}
\maketitle
\begin{abstract}
In Cator and Lopuha\"a~\cite{catorlopuhaa2009} an asymptotic expansion for the MCD estimators
is established in a very general framework.
This expansion requires the existence and non-singularity of the derivative in a first-order Taylor expansion.
In this paper, we prove the existence of this derivative for multivariate distributions that have a density
and provide an explicit expression.
Moreover, under suitable symmetry conditions on the density, we show that this derivative is non-singular.
These symmetry conditions include the elliptically contoured multivariate location-scatter model,
in which case we show that the minimum covariance determinant (MCD) estimators of multivariate location and covariance
are asymptotically equivalent to a sum of independent identically distributed vector and matrix valued
random elements, respectively.
This provides a proof of asymptotic normality and a precise description of the limiting covariance structure
for the MCD estimators.
\end{abstract}

\section{Introduction}
\label{sec:intro}
The MCD estimator~\cite{rousseeuw85} is one of the most popular robust methods to estimate multivariate location and scatter parameters.
These estimators, in particular the covariance estimator, also serve as robust plug-ins in other multivariate statistical techniques,
such as
principal component analysis~\cite{crouxhaesbroeck2000,serneelsverdonck2008},
multivariate linear regression~\cite{agullocrouxvanaelst2008,rousseeuwvanaelstvandriessenagullo2004},
discriminant analysis~\cite{hawkinsmclachlan1997},
factor analysis~\cite{pisonrousseeuwfilzmosercroux2003},
canonical correlations~\cite{taskinencrouxkankainen2006,zhou2009} and
error-in-variables models~\cite{fekriruizgazen2004},
among others (see also \cite{hubertrousseeuwvanaelst} for a more extensive overview).
For this reason, the distributional and the robustness properties of the MCD estimators
are essential for conducting inference and perform robust estimation in several statistical models.

The MCD estimators have the same high breakdown point as the minimum volume ellipsoid estimators
(e.g., see \cite{agullocrouxvanaelst2008,lopuhaarousseeuw1991}).
The asymptotic properties have first been studied by Butler, Davies and Jhun~\cite{butlerdaviesjuhn93}
in the framework of unimodal elliptically contoured densities,
who showed that the MCD location estimator converges at $\sqrt{n}$-rate towards
a normal distribution with mean equal to the MCD location functional.
In the same framework, Croux and Haesbroeck~\cite{crouxhaesbroeck99} give the expression for the influence function
of the MCD covariance functional and use this to compute limiting variances of the MCD covariance estimator.
The asymptotic theory has been extended and generalized in Cator and Lopuha\"a~\cite{catorlopuhaa2009},
who studied the MCD estimators and the corresponding functional in a very general framework.
They establish an asymptotic expansion of the type
\begin{equation}
\label{eq:expansion estimators}
\widehat{\theta}_n-\theta_0
=-\Lambda'(\theta_0)^{-1}\frac1n\sum_{i=1}^n\left(\Psi(X_i,\theta_0)-\mathds{E}\Psi(X_i,\theta_0)\right)+o_\mathbb{P}(n^{-1/2}),
\end{equation}
where $\widehat{\theta}_n$ and $\theta_0$ denote vectors consisting of the MCD estimators and
the MCD functional at the underlying distribution, respectively,
and $\Psi(\cdot,\theta_0)$ is a function that we will specify later on.
In principle, from this expansion a central limit theorem for the MCD estimator can be derived.
However, the expansion requires the existence and non-singularity of $\Lambda'(\theta_0)$.
Moreover, a more explicit expression of its inverse is desirable from a practical point of view,
since it determines the limiting variances.

In this paper we show that $\Lambda'(\theta_0)$ exists as long as the underlying distribution $P$ has a density $f$.
Its expression given in Theorem~\ref{th:derivative} offers the possibility to estimate the limiting variances of the MCD estimators
in any model where $P$ has a density.
We will also provide sufficient symmetry conditions on $f$ for $\Lambda'(\theta_0)$ to be non-singular.
This includes the special case of elliptically contoured densities
\[
f(x)=\det(\Sigma)^{-1/2}h((x-\mu)\Sigma^{-1}(x-\mu)),
\]
for which we show that the MCD location and covariance estimator are asymptotically equivalent to a sum of independent
vector and matrix valued random elements, respectively.
This exact expansion shows that at elliptically contoured densities the MCD location and MCD covariance estimator are
asymptotically independent and yields an explicit central limit theorem for both MCD estimators separately,
in such a way that the limiting covariances between elements of the location and covariance estimators
can be obtained directly from the covariances between elements of the summands.
Furthermore, the expansion for the MCD estimators is needed to obtain the limiting distribution of robustly re-weighted least squares estimators
for $(\mu,\Sigma)$, if  one uses the MCD estimators to assign the weights (see~\cite{lopuhaa1999}).

The paper is organized as follows.
In Section~\ref{sec:def} we define the MCD estimators and MCD functionals and discuss some results from~\cite{catorlopuhaa2009}
that are relevant for our setup.
In Section~\ref{sec:derivative} we establish the expression for the $\Lambda'(\theta_0)$ in terms of a linear mapping
and show that this mapping is non-singular under suitable symmetry conditions.
The special case of elliptically contoured densities is considered in Section~\ref{sec:elliptical contoured},
where we obtain an explicit expression of $\Lambda'(\theta_0)^{-1}$.
From this we derive an asymptotic expansion for the estimators, prove asymptotic normality,
and derive the influence function of the MCD functionals.
As special cases we recover results from~\cite{butlerdaviesjuhn93} and~\cite{crouxhaesbroeck99}
under weaker conditions.

All proofs have been postponed to an appendix at the end of the paper.

\section{Definition and preliminaries}
\label{sec:def}
For a sample $X_1,X_2,\ldots,X_n$ from a distribution $P$ on $\R^k$
the MCD estimator is defined as follows.
Fix a fraction $0<\gamma\leq 1$ and consider subsamples $S\subset \{X_{1},\ldots,X_{n}\}$
that contain $h_n\geq \lceil{n\gamma}\rceil$ points.
Define a corresponding trimmed sample mean and sample covariance matrix by
\begin{equation}
\label{eq:def MCD estimator}
\begin{split}
\widehat{T}_n(S) &=
\frac{1}{h_n}\sum_{X_i\in S}X_i,
\\
\widehat{C}_n(S)
&=
\frac{1}{h_n}\sum_{X_i\in S}(X_i-\widehat{T}_n(S))(X_i-\widehat{T}_n(S))'.
\end{split}
\end{equation}
Note that each subsample $S$ determines an ellipsoid
$E(\widehat{T}_n(S),\widehat{C}_n(S),\widehat{r}_n(S))$,
where for each $\mu\in\R^k$, $\Sigma$ symmetric positive definite, and $\rho>0$,
\begin{equation}
\label{eq:def ellipsoid}
E(\mu,\Sigma,\rho)
=
\big\{
x\in\R^k:(x-\mu)'\Sigma^{-1}(x-\mu)\leq\rho^2
\big\},
\end{equation}
and
\begin{equation}
\label{eq:def r_n}
\widehat{r}_n(S) = \inf\left\{ s>0\ : P_n\left(E(\widehat{T}_n(S),\widehat{C}_n(S),s)\right)\geq \gamma\right\},
\end{equation}
where $P_n$ denotes the empirical measure corresponding to the sample.
Let $S_n$ be a subsample that minimizes $\det(\widehat{C}_n(S))$
over all subsamples of size $h_n\geq \lceil{n\gamma}\rceil$, then the pair $(\widehat{T}_n(S_n),\widehat{C}_n(S_n))$
is an MCD-estimator.
Note that a minimizing subsample always exists, but it need not be unique.
In~\cite{catorlopuhaa2009} it is shown that a minimizing subsample $S_n$ always has exactly $\lceil n\gamma\rceil$ points
and is contained in the ellipsoid $E(\widehat{T}_n(S_n),\widehat{C}_n(S_n),\widehat{r}_n(S_n))$,
which separates $S_n$ from all other points in the sample.
Note that in~\cite{butlerdaviesjuhn93} (among others) one minimizes over subsamples of size $\lfloor n\gamma\rfloor$.
This is somewhat unnatural, since it may lead to subsamples $S$ for which $P_n(S)<\gamma$.
Moreover, it may lead to situations where the trimmed subsample does not contain the majority of the points,
e.g., if $\gamma=1/2$ and $n$ is odd, then $\lfloor n\gamma\rfloor=(n-1)/2$.
By considering subsamples~$S$ of size $h_n\geq \lceil n\gamma\rceil$ in definition~\eqref{eq:def MCD estimator},
we always have $P_n(S)\geq \gamma$ and for any $1/2\leq \gamma\leq 1$, the subsample contains the majority of points.

We define the MCD functionals in similar fashion.
Define a trimmed mean and covariance as follows,
\begin{equation}
\label{eq:def functionals}
\begin{split}
T_P(\phi)
&=
\frac1{\int \phi\,dP}\int x \phi(x)\,P(dx),\\
C_P(\phi)
&= \frac1{\int \phi\,dP}\int (x-T_P(\phi))(x-T_P(\phi))'\phi(x)\,P(dx).
\end{split}
\end{equation}
and define
\begin{equation*}
\label{eq:def r_P}
r_P(\phi) = \inf\left\{ s>0\ : P\left(E(T_P(\phi),C_P(\phi),s)\right)\geq \gamma\right\}.
\end{equation*}
for measurable $\phi:\R^k\to[0,1]$, such that $\int\phi\,dP\geq \gamma$ and $\int \|x\|^2\phi(x)\,P(dx)<\infty$.
Note that for $P=P_n$ and $\phi = \mathds{1}_{S}$,
we recover~\eqref{eq:def MCD estimator} and \eqref{eq:def r_n}.
If $\phi_P$ minimizes $\det(C_P(\phi))$ over all~$\phi$ considered above, then the pair $(T_P(\phi_P),C_P(\phi_P))$ is called an MCD functional.
In \cite{catorlopuhaa2009} it is shown that such a $\phi_P$ always exists and a characterization of a minimizing $\phi$ is provided.
From this characterization (Theorem~3.2 in~\cite{catorlopuhaa2009}) it follows that if $P$ has a density,
then
\begin{equation}\label{eq:phi=ellips}
\phi_P = \mathds{1}_{E_P}\quad \mbox{and}\quad P\left(E_P\right) = \gamma,
\end{equation}
where $E_P=E(T_P(\phi_P),C_P(\phi_P),r_P(\phi_P))$.
This means that the MCD functional defined by \eqref{eq:def functionals} coincides with the definition
through minimization over bounded Borel sets given in~\cite{butlerdaviesjuhn93}.

Throughout the paper we will assume that the MCD functional at $P$ is uniquely defined,
and we write $(\mu_0,\Sigma_0)=(T_P(\phi_P),C_P(\phi_P))$ and $\rho_0=r_P(\phi_P)$.
This holds, for instance, if $P$ has a unimodal elliptically contoured density (see Theorem~1 in~\cite{butlerdaviesjuhn93}).
We will also assume that $P$ has a density $f$ that satisfies the following condition:
\begin{quote}
\begin{itemize}
\item[(B)]
$f$ is continuous and strictly positive on a open neighborhood of the boundary of $E(\mu_0,\Sigma_0,\rho_0)$.
\end{itemize}
\end{quote}
In that case, it follows from Theorem~4.2 in~\cite{catorlopuhaa2009} that $\widehat{\theta}_n\to\theta_0$ with probability one,
where
\begin{equation}
\begin{split}
\label{eq:def theta}
\widehat{\theta}_n
&=
\big(\widehat{T}_n(S_n),\widehat{B}_n(S_n),\widehat{r}_n(S_n)\big),
\quad\text{with }
\widehat{B}_n(S_n)^2=\widehat{C}_n(S_n),
\\
\theta_0
&=
\big(\mu_0,\Gamma_0,\rho_0\big),
\quad\text{with }
\Gamma_0^2=\Sigma_0.
\end{split}
\end{equation}
Moreover, Theorem~5.1 in~\cite{catorlopuhaa2009} implies that expansion~\eqref{eq:expansion estimators} holds,
where $\Psi=(\Psi_1,\Psi_2,\Psi_3)$, is defined as
\begin{equation}
\label{eq:def Psi}
\begin{split}
\Psi_1(y,\theta)&= \mathds{1}_{\{\|G^{-1}(y-m)\|\leq r\}}G^{-1}(y-m)\\
\Psi_2(y,\theta)&= \mathds{1}_{\{\|G^{-1}(y-m)\|\leq r\}}\Big(G^{-1}(y-m)(y-m)'G^{-1}-I_k\Big)\\
\Psi_3(y,\theta)&= \mathds{1}_{\{\|G^{-1}(y-m)\|\leq r\}}-\gamma,
\end{split}
\end{equation}
and $\Lambda=(\Lambda_1,\Lambda_2,\Lambda_3)$, with
\begin{equation}
\label{eq:def Lambda}
\Lambda_j(\theta)=\int \Psi_j(y,\theta)\,P(dy),
\quad
\text{for }j=1,2,3,
\end{equation}
for $\theta=(m,G,r)$, with $y,t\in\R^k$, $r>0$, and $G\in \mathrm{PDS}(k)$.
Here, $\mathrm{PDS}(k)$ denotes the space of all positive definite symmetric $k\times k$ matrices.

\section{Existence and non-singularity of $\Lambda'(\theta_0)$}
\label{sec:derivative}
Let $\theta_0=(\mu_0,\Gamma_0,\rho_0)$ be the MCD functional at $P$.
Due to the indicator function in the expression of $\Psi(x,\theta_0)$ it can be seen, that the existence
of a derivative of $\Lambda(\theta)$ at $\theta_0$ cannot be expected in general if $P$
does not satisfy condition (B).
If $P$ does satisfy (B), then the derivative will depend on the behavior of $f$ on the boundary
of $E(\mu_0,\Sigma_0,\rho_0)$.
For $\rho>0$ and $\mu\in\R^k$, define
\begin{equation*}
\label{eq:def B(mu,rho)}
B(\mu,\rho)
=
\left\{x\in \R^k : \|x-\mu\|\leq \rho\right\},
\end{equation*}
and let $\sigma_0$ denote the surface measure on the boundary~$\partial B(0,\rho_0)$.
To describe the derivative of $\Lambda(\theta)$ at $\theta_0$, we introduce the measure
\begin{equation}
\label{eq:def nu}
\nu(d\omega)=\det(\Gamma_0)f(\Gamma_0\omega+\mu_0)\sigma_0(d\omega)
\quad
\text{for }\omega\in\partial B(0,\rho_0).
\end{equation}
Here, $\sigma_0$ denotes the Lebesgue surface measure on $\partial B(0,\rho_0)$.

\bigskip

Note that our parameter $\theta_0$ is an element of $\Theta = \R^k\times {\rm PDS}(k)\times \R$. This means that the derivative of $\Lambda$ at $\theta_0$, if it exists, can be described as a linear mapping on the tangent space of $\Theta$ in $\theta_0$, which is given by $\R^k\times S(k)\times \R$.
Here, $S(k)$ denotes the space of all symmetric $k\times k$ matrices.
The derivatives of $\Lambda_1$, $\Lambda_2$ and $\Lambda_3$ are given as linear mappings by the following theorem.
\begin{theorem}
\label{th:derivative}
Suppose that $P$ satisfies (B)
and let the MCD functional $\theta_0=(\mu_0,\Gamma_0,\rho_0)$ be uniquely defined at $P$.
For $j=1,2,3$,
the derivatives of $\Lambda_j$ are given by the following linear mappings, with $(h,A,s)\in \R^k\times S(k)\times \R$:
\[
\begin{split}
\Lambda_1'(\theta_0)(h,A,s)
&=
-\gamma \Gamma_0^{-1}h
+
\int_{\partial B_0}
\left(
\frac{\omega'\Gamma_0^{-1}h}{\rho_0} +\frac{\omega'(\Gamma_0^{-1}A+A\Gamma_0^{-1})\omega}{2\rho_0}+s
\right)
\omega\, \nu(d\omega)\\
\Lambda_2'(\theta_0)(h,A,s)
&=
-\gamma(\Gamma_0^{-1}A+A\Gamma_0^{-1})
+
\int_{\partial B_0}
\left(
\frac{\omega'\Gamma_0^{-1}h}{\rho_0} +\frac{\omega'(\Gamma_0^{-1}A+A\Gamma_0^{-1})\omega}{2\rho_0}+s
\right)
\Big(\omega\omega'-I\Big)\, \nu(d\omega)\\
\Lambda_3'(\theta_0)(h,A,s)
&=
\int_{\partial B_0}
\left(
\frac{\omega'\Gamma_0^{-1}h}{\rho_0} +\frac{\omega'(\Gamma_0^{-1}A+A\Gamma_0^{-1})\omega}{2\rho_0}+s
\right)\,
\nu(d\omega),
\end{split}
\]
where $B_0=B(0,\rho_0)$ and $\nu(d\omega)$ is defined in~\eqref{eq:def nu}.
\end{theorem}
Note that Theorem~\ref{th:derivative} also has practical implications.
According to Theorem~5.1 in~\cite{catorlopuhaa2009}, the MCD estimator $\widehat{\theta}_n=(\widehat{\mu}_n,\widehat{\Gamma}_n,\widehat{\rho}_n)$,
represented as a vector, is asymptotically normal with mean zero and limiting variance given by the
covariance matrix of $Z=\Lambda'(\theta_0)^{-1}\Psi(X_1,\theta_0)$.
This means that the  expression for $\Lambda'(\theta_0)$ enables one to estimate the limiting variance of the MCD estimators
in any model where $P$ has a density,
which goes far beyond the traditional elliptically contoured densities.
An estimate for~$\Lambda'(\theta_0)$ can obtained by plugging in the estimate $\widehat{\theta}_n$ for $\theta_0$,
replacing $\sigma_0$ by the surface measure $\widehat{\sigma}_n$ on $\partial B(0,\widehat{\rho}_n)$, and
using a nonparametric estimate for the density $f$ on the boundary of $B(0,\widehat{\rho}_n)$,
e.g., histogram or kernel type estimates.
As long as the estimate $\widehat{\Lambda}'(\widehat{\theta}_n)$ turns out to be non-singular, the limiting covariance matrix of
$\sqrt{n}(\widehat{\theta}_n-\theta_0)$ can be estimated by the sample covariance of
the $Z_i=\widehat{\Lambda}'(\widehat{\theta}_n)^{-1}\Psi(X_i,\widehat{\theta}_n)$.

\bigskip

We proceed by finding sufficient conditions for $\Lambda'(\theta_0)$ to be non-singular.
We would have non-singularity, if for all $(h,A,s)\in \R^k\times S(k)\times \R$,  $\Lambda'(\theta_0)(h,A,s)=0$ implies that $(h,A,s)=(0,0,0)$.
From the expressions in Theorem~\ref{th:derivative} it can be seen that this cannot be expected without
further assumptions on $f$.
Suitable symmetry assumptions on $f$ will simplify the expressions for the derivative, in
which case non-singularity can be established.
First note that without any further assumptions on~$f$,
we do always have the following result.
\begin{lemma}
\label{lem:trace=0}
Suppose that $P$ satisfies (B)
and let the MCD functional $\theta_0=(\mu_0,\Gamma_0,\rho_0)$ be uniquely defined at $P$.
Let $\Lambda$ be defined by \eqref{eq:def Lambda} and suppose $\Lambda'(\theta_0)(h,A,s)=0$.
Then $\mathrm{Tr}(\Gamma_0^{-1}A)=0$.
\end{lemma}
Next, we consider the case where the density $f$ is point symmetric with respect to the center of $E(\mu_0,\Sigma_0,\rho_0)$,
i.e.,
\begin{equation}
\label{eq:point symmetry}
f(-\Gamma_0\omega+\mu_0)=f(\Gamma_0\omega+\mu_0),
\quad
\text{for }
\omega\in\partial B(0,\rho_0).
\end{equation}
In that case we have the following lemma.
\begin{lemma}
\label{lem:point symmetry}
Suppose $P$ satisfies (B) and \eqref{eq:point symmetry},
and let the MCD functional $\theta_0=(\mu_0,\Gamma_0,\rho_0)$ be uniquely defined at $P$.
Let $\Lambda$ be defined by \eqref{eq:def Lambda} and suppose that $\Lambda'(\theta_0)(h,A,s)=0$.
Then
\begin{equation}
\label{eq:relation s-A}
s=-\frac{1}{2\rho_0\nu_0}\int_{\partial B_0}
\omega'(\Gamma_0^{-1}A+A\Gamma_0^{-1})\omega\,\nu(d\omega),
\end{equation}
where $B_0=B(0,\rho_0)$ and $\nu_0=\nu(\partial B_0)$, with $\nu$ defined in~\eqref{eq:def nu}.
If, in addition, for all $i=1,2\ldots,k$
\begin{equation}
\label{eq:cond h=0}
\int_{\partial B_0}\omega_i^2\,\nu(d\omega)\ne \gamma\rho_0,
\end{equation}
then $h=0$.
\end{lemma}
Point symmetry will not be sufficient to conclude that $A=0$ from $\Lambda'(\theta_0)(h,A,s)=0$.
The slightly stronger condition of half-space symmetry will suffice, i.e.,
\begin{equation}
\label{eq:half space symmetry}
f(\Gamma_0\omega_{(-i)}+\mu_0)=f(\Gamma_0\omega+\mu_0),
\quad
\text{ where }
\omega_{(-i)}=(\omega_1,\ldots,\omega_{i-1},-\omega_i,\omega_{i+1},\ldots,\omega_k),
\end{equation}
for all $i=1,2,\ldots,k$
and $\omega\in\partial B(0,\rho_0)$.
To describe sufficient conditions for non-singularity, we define
the matrix $M$ with elements
\begin{equation}
\label{eq:def M}
M_{ij}
=
\int_{\partial B(0,\rho_0)}\omega_i^2\omega_{j}^2\,\nu(d\omega)
-
\frac{1}{\nu_0}
\int_{\partial B(0,\rho_0)}\omega_i^2\,\nu(d\omega)
\int_{\partial B(0,\rho_0)}\omega_j^2\,\nu(d\omega)
-
2\gamma \rho_0\mathds{1}_{\{i=j\}},
\end{equation}
for $i,j=1,2,\ldots,k$,
where $\nu_0=\nu(\partial B(0,\rho_0))$ and $\nu(d\omega)$ is defined by \eqref{eq:def nu}.
We then have the following theorem.
\begin{theorem}
\label{th:non-singular}
Suppose that $P$ satisfies (B) and \eqref{eq:half space symmetry},
and let the MCD functional $\theta_0=(\mu_0,\Gamma_0,\rho_0)$ be uniquely defined at $P$.
Suppose that \eqref{eq:cond h=0} holds, that for all $i,j=1,2\ldots,k$ with $i\ne j$,
\begin{equation}
\label{eq:cond Anm=0}
\int_{\partial B(0,\rho_0)}\omega_i^2\omega_j^2\,\nu(d\omega)
\ne
\gamma \rho_0,
\end{equation}
where $\nu$ is defined in~\eqref{eq:def nu}, and that the matrix $M$ defined in \eqref{eq:def M} is such that
for any $x\in\R^k$,
\begin{equation}
\label{eq:property M}
Mx=0
\text{ and }
x_1+\cdots+x_k=0
\quad\Rightarrow\quad
x=0.
\end{equation}
Then, for $\theta_0=(\mu_0,\Gamma_0,\rho_0)$, the derivative $\Lambda'(\theta_0)$ is non-singular as a linear map on $\R^k\times S(k)\times \R$.
\end{theorem}
Example of densities that satisfy \eqref{eq:half space symmetry} are elliptically contoured densities.
However, also affine transformations of densities that have independent marginal densities that are symmetric around zero,
i.e.,
\[
f(x)=g(\Gamma^{-1}(x-\mu)),
\quad
\text{where }
g(x_1,\ldots,x_k)=g_1(x_1)\cdots g_k(x_k)
\text{ and }
g_i(x_i)=g_i(-x_i),
\]
satisfy \eqref{eq:half space symmetry}.

\section{Elliptically contoured densities}
\label{sec:elliptical contoured}
Suppose that $P$ has an elliptically contoured density,
i.e.,
\begin{equation}
\label{eq:def ellip density}
f(x)=\det(\Sigma)^{-1/2}h((x-\mu)\Sigma^{-1}(x-\mu))
\quad
\text{where }\mu\in\R^k, \Sigma\in \mathrm{PDS}(k),
\end{equation}
and $h:[0,\infty)\to[0,\infty)$ is decreasing so that $P$ is unimodal.
In this case, it follows from the characterization for the $\phi$ function that minimizes $\det(C_P(\phi))$
(see Theorem~3.2 in~\cite{catorlopuhaa2009}),
that our definition of the MCD functional coincides with the one in~\cite{butlerdaviesjuhn93},
who show that the MCD functionals are unique:
\begin{equation}
\label{eq:MCD at ellip cont}
\mu_0=\mu,
\quad
\Sigma_0=\alpha(\gamma)^2\Sigma,
\quad\text{and}\quad
\rho_0^2=\frac{r(\gamma)^2}{\alpha(\gamma)^2},
\end{equation}
where
\begin{equation}
\label{eq:def alpha}
\alpha(\gamma)^2
=
\frac{2\pi^{k/2}}{\gamma k\Gamma(k/2)}\int_0^{r(\gamma)}h(r^2)r^{k+1}\,dr,
\end{equation}
and where $r(\gamma)$ is determined by
\begin{equation}
\label{eq:def r(gamma)}
\frac{2\pi^{k/2}}{\Gamma(k/2)}\int_0^{r(\gamma)}h(r^2)r^{k-1}\,dr=\gamma.
\end{equation}
The next proposition shows that for elliptically contoured densities the derivative $\Lambda'(\theta_0)$
exists and is non-singular.
\begin{prop}
\label{prop:cond elliptical}
Let $P$ have an elliptically contoured density as defined in~\eqref{eq:def ellip density}
with $h$ non-increasing such that $P$ is unimodal.
Then all conditions of Theorem~\ref{th:non-singular} are satisfied.
\end{prop}
We proceed by obtaining asymptotic expansions for the MCD estimators in the case of
elliptically contoured densities.
Because the estimators are affine equivariant, it suffices to consider the
spherically symmetric case $(\mu,\Sigma)=(0,I)$.
The next theorem provides the expressions for~$\Lambda'(\theta_0)$ and its inverse at spherically symmetric densities.
\begin{theorem}
\label{th:derivative ellip}
Let $P$ have an spherically symmetric density $f(x)=h(\|x\|^2)$ with $h$ decreasing such that $P$ is unimodal.
Let $r=r(\gamma)$ and $\alpha=\alpha(\gamma)$ be defined in~\eqref{eq:def r(gamma)} and~\eqref{eq:def alpha}, respectively,
and let $D=\Lambda'(\theta_0)$, for $\theta_0=(\mu_0,\Gamma_0,\rho_0)=(0,\alpha I,r/\alpha)$.
Then the linear mapping~$D$ is given by
\[
\begin{split}
D_1:(h,A,s)
&\mapsto
\beta_1h\\
D_2:(h,A,s)
&\mapsto
\beta_2 A+ \beta_3 \mathrm{Tr}(A)\cdot I+ \beta_4s\cdot I\\
D_3:(h,A,s)
&\mapsto
\beta_5\mathrm{Tr}(A)+\beta_6s,
\end{split}
\]
and the inverse linear mapping $D^{\mathrm{inv}}$ is given by
\[
\begin{split}
\left[D^{\mathrm{inv}}\right]_1:(g,B,t)
&\mapsto
\beta_1^{-1}g\\
\left[D^{\mathrm{inv}}\right]_2:(g,B,t)
&\mapsto
\beta_2^{-1}B
+\frac{\alpha(\beta_3\beta_6-\beta_4\beta_5)}{2\gamma\beta_2\beta_6}\mathrm{Tr}(B)\cdot I
+\frac{\alpha\beta_4}{2\gamma\beta_6}t\cdot I\\
\left[D^{\mathrm{inv}}\right]_3:(g,B,t)
&\mapsto
\frac{\alpha\beta_5}{2\gamma\beta_6}\mathrm{Tr}(B)-\frac{\alpha(\beta_2+k\beta_3)}{2\gamma\beta_6}t
\end{split}
\]
where
\[
\begin{array}{ll}
\beta_1
=
\dfrac1{\alpha}
\left(
\dfrac{\rho_0}{k}\nu_0-\gamma
\right)<0,
&\quad
\beta_4
=
\dfrac{\rho_0^2}{k}\nu_0
-
\nu_0,\\
\\
\beta_2
=
\dfrac{2\rho_0^3\nu_0}{\alpha k(k+2)}
-
\dfrac{2\gamma}{\alpha}<0,
&\quad
\beta_5
=
\dfrac{\rho_0\nu_0}{k\alpha },\\
\\
\beta_3
=
\dfrac{\rho_0^3\nu_0}{\alpha k(k+2)}
-
\dfrac{\rho_0\nu_0}{k\alpha},
&\quad
\beta_6
=
\nu_0>0,
\end{array}
\]
with $B_0=B(0,\rho_0)$ and
\[
\nu_0
=
\nu(\partial B_0)
=
\frac{2\pi^{k/2}}{\Gamma(k/2)}h(r^2) r^{k-1}\alpha.
\]
\end{theorem}
An immediate consequence of Theorem~\ref{th:derivative ellip} is the next corollary,
which shows that the MCD estimators of location and covariance are asymptotically equivalent
to a sum of independent identically distributed vector and matrix valued random elements, respectively.
\begin{corollary}
\label{cor:expansion}
Suppose that $P$ has a spherically symmetric density $f(x)=h(\|x\|^2)$ with $h$ decreasing such that $P$ is unimodal.
Let $r=r(\gamma)$ and $\alpha=\alpha(\gamma)$ be defined in~\eqref{eq:def r(gamma)} and~\eqref{eq:def alpha}, respectively.
Then for $n\to\infty$,
\[
\begin{split}
\sqrt{n}\widehat{\mu}_n
&=
\frac1{\sqrt{n}}\sum_{i=1}^n
\pi\mathds{1}_{\{\|X_i\|\leq r\}}X_i+o_\mathds{P}(1);\\
\sqrt{n}
(\widehat{\Sigma}_n-\alpha^2 I)
&=
\frac{1}{\sqrt{n}}
\sum_{i=1}^n
\left[
\mathds{1}_{\{\|X_i\|\leq r\}}
\left(
\kappa_1\cdot I
+
\kappa_2
\|X_i\|^2\cdot I
+
\kappa_3
X_iX_i'
\right)
+
\kappa_4\cdot I
\right]
+o_\mathbb{P}(1);\\
\sqrt{n}
\left(
\widehat{\rho}_n-\frac{r}{\alpha}
\right)
&=
\frac{1}{\sqrt{n}}
\sum_{i=1}^n
\left[
\lambda_1\mathds{1}_{\{\|X_i\|\leq r\}}\|X_i\|^2
+
\lambda_2\mathds{1}_{\{\|X_i\|\leq r\}}
+
\lambda_3
\right]
+o_\mathbb{P}(1),
\end{split}
\]
where $\pi=-(\alpha\beta_1)^{-1}$ and
\[
\begin{split}
\kappa_1
&=
-\dfrac{r^2}{k\gamma},
\quad
\kappa_2
=
\dfrac{\alpha\beta_2+2\gamma}{k\gamma\alpha\beta_2},
\quad
\kappa_3
=
-\dfrac{2}{\alpha\beta_2},
\quad
\kappa_4
=
\dfrac{r^2-k\alpha^2}{k}
\\
\lambda_1
&=
-\frac{r}{2k\gamma\alpha^3},
\quad
\lambda_2=\frac{r^3}{2k\gamma\alpha^3}-\frac{1}{\beta_6},
\quad
\lambda_3=\frac{\gamma}{\beta_6}+\frac{r}{2k\alpha^3}\left(k\alpha^2-r^2\right),
\end{split}
\]
with $\beta_1$, $\beta_2$ and $\beta_6$ defined in Theorem~\ref{th:derivative ellip}.
\end{corollary}
We proceed by obtaining the limit distribution of the MCD estimators.
To describe the limiting distribution of a random matrix
we use the operator $\text{vec}(\cdot)$ which
stacks the columns of a matrix $M$ on top of each other, i.e.,
\[
\text{vec}(M) =(M_{11},\ldots,M_{1k},\ldots,M_{k1},\ldots,M_{kk})'.
\]
We will also need the commutation matrix $C_{k,k}$,
which is a $k^2\times k^2$ matrix consisting of $k\times k$ blocks:
$C_{k,k}=(\Delta_{ij})_{i,j=1}^k$,
where each $(i,j)$-th block is equal to a $k\times k$-matrix
$\Delta_{ji}$, which is~1 at entry $(j,i)$ and 0 everywhere else.
Finally, for matrices $M$ and $N$, the Kronecker product $M\otimes N$ is a $k^2\times k^2$
matrix consisting of $k\times k$ blocks, with the $(i,j)$-th block equal to $m_{ij}N$.
\begin{theorem}
\label{th:AN}
Suppose that $P$ has a spherically symmetric density $f(x)=h(\|x\|^2)$ with $h$ decreasing such that $P$ is unimodal.
Let $r=r(\gamma)$ and $\alpha=\alpha(\gamma)$ be defined in~\eqref{eq:def r(gamma)} and~\eqref{eq:def alpha}, respectively.
Let $\widehat{\mu}_n$, $\widehat{\Sigma}_n$ and $\widehat{\rho}_n$ the MCD estimators.
Then
\begin{itemize}
\item[(i)]
$\widehat{\mu}_n$ and $(\widehat{\Sigma}_n,\widehat{\rho}_n)$ are asymptotically independent,
the diagonal elements of
$\widehat{\Sigma}_n$ are asymptotically independent from the off-diagonal elements and $\widehat{\rho}_n$,
and the off-diagonal elements of $\widehat{\Sigma}_n$ are asymptotically mutually independent;
\item[(ii)]
$\sqrt{n}\widehat{\mu}_n$ is asymptotically normal with mean zero and covariance matrix $\tau I$,
where
\[
\tau=\frac{k^2\gamma\alpha^4}{(k\gamma\alpha-r\nu_0)^2},
\]
where $\nu_0$ is defined in Theorem~\ref{th:derivative ellip};
\item[(iii)]
$\sqrt{n}(\mathrm{vec}(\widehat{\Sigma}_n)-\alpha^2\mathrm{vec}(I))$ is asymptotically normal with mean zero and covariance matrix
\[
\sigma_1(I+C_{k,k})(I\otimes I)+\sigma_2\mathrm{vec}(I)\mathrm{vec}(I)',
\]
where
\[
\begin{split}
\sigma_1&=\dfrac{\kappa_3^2}{k(k+2)}\mathbb{E}\mathds{1}_{\{\|X_1\|\leq r\}}\|X_1\|^4\\
\sigma_2
&=
-\frac2k\sigma_1
+
\frac{1}{k^2\gamma^2}
\mathbb{E}\mathds{1}_{\{\|X_1\|\leq r\}}\|X_1\|^4
-
\frac{\gamma r^4-2k\gamma r^2\alpha^2+k^2\gamma\alpha^4+2kr^2\alpha^2-r^4}{\gamma k^2}
\end{split}\]
where $\kappa_3$ is defined in Corollary~\ref{cor:IF};
\item[(iv)]
$\sqrt{n}(\widehat{\rho}_n-r/\alpha)$ is asymptotically normal with mean zero and variance
\[
\sigma_\rho^2=
\lambda_1^2\mathbb{E}\mathds{1}_{\{\|X_1\|\leq r\}}\|X_1\|^4
+
\frac{2k^2\nu_0r\alpha^5\gamma-\nu_0^2kr^4\alpha^2+4k^2\gamma^2\alpha^6-4kr^3\nu_0\alpha^3\gamma+r^6\nu_0^2}{4k^2\alpha^6\nu_0^2\gamma },
\]
where $\lambda_1$ is defined in Corollary~\ref{cor:expansion}
and $\nu_0$ is defined in Theorem~\ref{th:derivative ellip}.
\end{itemize}
\end{theorem}
With Theorem~\ref{th:AN}(i) we recover Theorem~4 in~\cite{butlerdaviesjuhn93}.
Note however, that the assumption of $h$ being differentiable (see~\cite{butlerdaviesjuhn93}) is not required in our approach.
Furthermore, it can be seen from the expression of the limiting variance of $\widehat{\Sigma}_n$ that
in the spherically symmetric case:
\[
\begin{split}
\sqrt{n}(\widehat{\Sigma}_{n,ii}-\alpha^2)&\to N(0,2\sigma_1+\sigma_2)\\
\sqrt{n}\widehat{\Sigma}_{n,ij}&\to N(0,\sigma_1)\\
\sqrt{n}
\left(
\begin{array}{c}
\widehat{\Sigma}_{n,ii}-\alpha^2
\\
\widehat{\Sigma}_{n,jj}-\alpha^2
\end{array}
\right)
&\to
N\left(
\left(
  \begin{array}{c}
    0 \\
    0 \\
  \end{array}
\right),
\left(
  \begin{array}{cc}
    2\sigma_1+\sigma_2 & \sigma_2 \\
    \sigma_2 & 2\sigma_1+\sigma_2 \\
  \end{array}
\right)
\right),
\qquad
i\ne j
\end{split}
\]
for $i,j=1,2,\ldots,k$.

Because $\widehat{\mu}_n$ and $\widehat{\Sigma}_n$ are affine equivariant the limiting
distributions for the MCD estimators in the case of general
$\mu\in\R^k$ and $\Sigma\in\mbox{PDS}(k)$ can be obtained easily.
Suppose that $X_1,\ldots,X_n$ are independent with density
$$
f(x)=\det(\Sigma)^{-1/2}h((x-\mu)^T\Sigma^{-1}(x-\mu)).
$$
Because of affine equivariance it follows immediately that
$\sqrt{n}(\widehat{\mu}_n-\mu)$ is asymptotically normal with zero
mean and covariance matrix
$\Gamma(\tau I)\Gamma=\tau\Sigma$,
where $\Gamma^2=\Sigma$.
Similarly $\sqrt{n}(\mathrm{vec}(\widehat{\Sigma}_n)-\alpha^2\mathrm{vec}(\Sigma))$ is asymptotically normal
with mean
zero and covariance matrix
$\mathbb{E}\mathrm{vec}(\Gamma M\Gamma )
\mathrm{vec}(\Gamma M\Gamma )'$,
where $M$ is the random matrix with
$\mathbb{E}\mathrm{vec}(M)\mathrm{vec}(M)'=
\sigma_1(I+C_{k,k})+\sigma_2\mathrm{vec}(I)\mathrm{vec})I)'$.
It follows from Lemma~5.2 in~\cite{lopuhaa89}, that
$$
\mathbb{E}\mathrm{vec}(\Gamma M\Gamma )
\mathrm{vec}(\Gamma M\Gamma )'
=
\sigma_1(I+C_{k,k})(\Sigma\otimes\Sigma)
+
\sigma_2\mathrm{vec}(\Sigma)\mathrm{vec}(\Sigma)'.
$$
This means that we have the following general corollary of
Theorem~\ref{th:AN}.
\begin{corollary}
\label{cor:asymp S general}
Suppose that $X_1,\ldots,X_n$ are independent with an elliptical
contoured density
$$
f(x)=\det(\Sigma)^{-1/2}h((x-\mu)'\Sigma^{-1/2}(x-\mu))
\quad,\mu\in\R^k,\,\Sigma\in\mathrm{PDS}(k),
$$
where $h:[0,\infty)\to[0,\infty)$ is non-increasing such that $f$ is unimodal.
Let $(\widehat{\mu}_n,\widehat{\Sigma}_n)$ be the MCD-estimators.
Then $\widehat{\mu}_n$ and $\widehat{\Sigma}_n$ are asymptotically independent,
$\sqrt{n}(\widehat{\mu}_n-\mu)$ has a limiting normal distribution with zero mean
and covariance matrix $\tau\Sigma$
and $\sqrt{n}(\mathrm{vec}(\widehat{\Sigma}_n)-\alpha^2\mathrm{vec}(\Sigma))$ has a limiting normal distribution with zero
mean and covariance matrix
$\sigma_1(I+C_{k,k})(\Sigma\otimes\Sigma)
+
\sigma_2\mathrm{vec}(\Sigma)\mathrm{vec}(\Sigma)'$,
where $\tau$, $\sigma_1$ and $\sigma_2$ are given in
Theorem~\ref{th:AN}.
\end{corollary}
Another corollary of Theorem~\ref{th:derivative ellip} is the expression for the influence function of the MCD functional.
The influence function of a functional~$\Theta(\cdot)$ at $P$ is defined as
\begin{equation}
\label{eq:def IF}
\mathrm{IF}(x,\Theta,P)
=
\lim_{\eps\downarrow 0}\frac{\Theta((1-\eps)P+\eps\delta_x)-\Theta(P)}{\eps},
\end{equation}
if this limit exists, where $\delta_x$ is the Dirac measure at $x\in\R^k$.
Denote by $\mu(P)=T_P(\phi_P)$, $\Sigma(P)=C_P(\phi_P)$, and~$\rho(P)=r_P(\phi_P)$ the MCD functionals at distribution $P$.
We then have the following corollary.
\begin{corollary}
\label{cor:IF}
Suppose that $P$ has a spherically symmetric density $f(x)=h(\|x\|^2)$ with $h$ decreasing such that $P$ is unimodal.
Let $r=r(\gamma)$ and $\alpha=\alpha(\gamma)$ be defined in~\eqref{eq:def r(gamma)} and~\eqref{eq:def alpha}, respectively.
Then, for $x\in\R^k$ such that $\|x\|\ne r$, the influence functions of the functionals $\mu(P)$, $\Sigma(P)$ and $\rho(P)$ are given by
\[
\begin{split}
\mathrm{IF}(x;\mu,P)
&=
\pi\mathds{1}_{\{\|x\|\leq r\}}x\\
\mathrm{IF}(x;\Sigma,P)
&=
\mathds{1}_{\{\|x\|\leq r\}}
\left(
\kappa_1\cdot I
+
\kappa_2
\|x\|^2\cdot I
+
\kappa_3
xx'
\right)
+
\kappa_4\cdot I
\\
\mathrm{IF}(x;\rho,P)
&=
\lambda_1\mathds{1}_{\{\|x\|\leq r\}}\|x\|^2
+
\lambda_2\mathds{1}_{\{\|x\|\leq r\}}
+
\lambda_3
\end{split}
\]
where $\pi$, $\kappa_1,\kappa_2,\kappa_3,\kappa_4$ and $\lambda_1,\lambda_2,\lambda_3$
are defined in Corollary~\ref{cor:expansion}.
\end{corollary}
Clearly, all the expressions in Corollary~\ref{cor:IF} are bounded uniformly for $\|x\|\ne r(\gamma)$.
For $x\in\R^k$ with $\|x\|=r(\gamma)$, it is not clear whether the limit in \eqref{eq:def IF} exists,
not even in the case of a unimodal spherically symmetric density.
As a special case of Corollary~\ref{cor:IF} we recover Theorem~1 in~\cite{crouxhaesbroeck99}.
However, we do not need the assumption that $h$ is differentiable (see~\cite{crouxhaesbroeck99}).
In order to see that our expressions coincide with the ones in~\cite{crouxhaesbroeck99},
note that their quantities $g$, $\alpha$, $q_\alpha$, and $c_\alpha$,
correspond to our $h$, $1-\gamma$, $r(\gamma)^2$, and $1/\alpha(\gamma)^2$, respectively,
and that they consider the Fisher consistent version of the covariance functional,
i.e., $c_\alpha\times \Sigma(P)$.
Moreover, their expression $b_1-kb_2$ is simply equal to 1.
For further discussion on $\mathrm{IF}(x;\Sigma,P)$ at spherically symmetric densities
and corresponding graphs,
we refer to~\cite{crouxhaesbroeck99}.

\section{Appendix}
\label{sec:appendix}
\subsection{Proofs for Section~\ref{sec:derivative}}
The proof of Theorem~\ref{th:derivative} requires the following lemma,
which helps to describe the derivative of~$\Lambda$ when $(\mu_0,\Gamma_0,s) = (0,I,r)$, in terms of a linear mapping.
Let $M(k)$ be the space of all $k\times k$ matrices.
\begin{lemma}
\label{lem:expansion boundary}
Let $r>0$ and $\phi:\R^k\to\R^m$, which is continuous on $\partial B(0,r)$.
Define the mapping $L:\R^k\times M(k)\times \R \to \R^m$ by
\[
L(h,A,s)=\int_{E(h,(I+A)(I+A)',r+s)}\phi(y)\,\mathrm{d}y.
\]
Then, the derivative of $L$ at $(h_0,A_0,s_0)=(0,0,0)$, is given by the continuous linear mapping
\[
L'(0,0,0)(h,A,s)
=
\int_{\partial B(0,r)}
\left(
\frac{\omega'h}{r}+\frac{\omega'(A+A')\omega}{2r}+s
\right)
\phi(\omega)\,\sigma_0(d\omega),
\]
with $(h,A,s)\in \R^k\times M(k)\times \R$.
\end{lemma}
\textbf{Proof:}
The derivative can be found as the sum of the derivatives of
\begin{equation}
\label{eq:three integrals}
\int_{B(h,r)}\phi(y)\,\mathrm{d}y,
\quad
\int_{E(0,(I+A)(I+A)',r)}\phi(y)\,\mathrm{d}y,
\quad\text{and }
\int_{B(0,r+s)}\phi(y)\,\mathrm{d}y.
\end{equation}
For the first integral, consider
\[
\int_{B(h,r)}\phi(y)\,\mathrm{d}y-\int_{B(0,r)}\phi(y)\,\mathrm{d}y
=
\int \left(\mathds{1}_{B(h,r)}-\mathds{1}_{B(0,r)}\right)\phi(y)\,\mathrm{d}y,
\]
for $\|h\|\to0$.
In first order this reduces to integration over $\partial B(0,r)$.
Let $\omega\in\partial B(0,r)$, let $v=(1+\delta)\omega\in \partial B(h,r)$,
and let $\alpha$ denote the angle between $\omega$ and $h$.
Then the law of cosines yields that
\[
r^2
=
\|v\|^2+\|h\|^2-2\|v\|\cdot\|h\|\cos\alpha
=
\left(
\|v\|-\frac{\omega'h}{\|\omega\|}
\right)^2
+
\|h\|^2(\sin\alpha)^2.
\]
Since $\|\omega\|=r$, in first order we find $r^2=(1+\delta)r^2-\omega'h$, or $\delta=(\omega'h)/r^2$.
This means that for each $\omega\in\partial B(0,r)$, the length over which we integrate $\phi(\omega)$ is
$\|v\|-\|\omega\|=\delta\|\omega\|=(\omega'h)/r$.
Since $\phi$ is continuous at $\omega\in\partial B(0,r)$, we get, for $\|h\|\to 0$,
\[
\int_{B(h,r)}\phi(y)\,\mathrm{d}y-\int_{B(0,r)}\phi(y)\,\mathrm{d}y
=
\int_{\partial B(0,r)} \frac{\omega'h}{r}\,\phi(\omega)\,\sigma_0(d\omega)+o(\|h\|).
\]
For the second integral in~\eqref{eq:three integrals} we consider
\[
\int \left(
\mathds{1}_{E(0,(I+A)(I+A)',r)}-\mathds{1}_{B(0,r)}\right)\,\phi(y)\,\mathrm{d}y,
\]
for $\|A\|\to0$,
which reduces to integration over $\omega\in\partial B(0,r)$.
Let $v=(1+\delta)\omega$ be such that $\|(I+A)^{-1}v\|=r$.
Then
\[
(1+\delta)^2=\frac{r^2}{\omega'(I+A')^{-1}(I+A)^{-1}\omega}.
\]
Since, for $\|A\|\to0$, we have $(I+A')^{-1}(I+A)^{-1}=I-A-A'+O(\|A\|^2)$,
it follows that
\[
\delta
=
\frac{\omega'(A+A')\omega}{2r^2}+O(\|A\|^2).
\]
This means that for each $\omega\in\partial B(0,r)$ the length over which we integrate $\phi(\omega)$ is
\[
\|v\|-\|\omega\|=\delta\|\omega\|=\frac{\omega'(A+A')\omega}{2r}+O(\|A\|^2).
\]
This implies that
\[
\int \left(
\mathds{1}_{E(0,(I+A)(I+A)',r)}-\mathds{1}_{B(0,r)}\right)\phi(y)\,\mathrm{d}y
=
\int_{\partial B(0,r)}
\frac{\omega'(A+A')\omega}{2r} \phi(\omega)\,\sigma_0(d\omega)+o(\|A\|).
\]
Finally, for the third integral in~\eqref{eq:three integrals} we obtain
\[
\int \left(
\mathds{1}_{B(0,r+s)}-\mathds{1}_{B(0,r)}\right)\phi(y)\,\mathrm{d}y
=
s\int_{\partial B(0,r)}\phi(\omega)\,\sigma_0(d\omega)+o(s).
\]
Summing the three linear mappings, yields the desired result.
\hfill\tqed

\bigskip

\noindent
\textbf{Proof of Theorem~\ref{th:derivative}:}
First note that everything can be rescaled to the situation with $\mu_0=0$ and $\Gamma_0=I$,
i.e., for any function~$g(y)$, we have
\begin{equation}
\label{eq:transformation}
\begin{split}
&
\int
\mathds{1}_{E(\mu_0+h,(\Gamma_0+A)(\Gamma_0+A)',\rho_0+s)}(y)\,g(y)\,P(dy)\\
&=
\det(\Gamma_0)
\int
\mathds{1}_{E(\widetilde{h},(I+\widetilde{A})^2,\rho_0+s)}(z)\,
g(\Gamma_0 z+\mu_0)f(\Gamma_0 z+\mu_0)\,\mathrm{d}z,
\end{split}
\end{equation}
where $\widetilde{h}=\Gamma_0^{-1}h$ and $\widetilde{A}=\Gamma_0^{-1}A$.
To compute $\Lambda_3'(\theta_0)$, take $g(y)=1$ in~\eqref{eq:transformation} and for $\eta=(h,A,s)\to (0,0,0)$, consider
\[
\begin{split}
\Lambda_3(\theta_0+\eta)-\Lambda_3(\theta_0)
&=
\int
\mathds{1}_{E(\mu_0+h,(\Gamma_0+A)(\Gamma_0+A)',\rho_0+s)}(y)
\,P(dy)
-
\int
\mathds{1}_{E(\mu_0,\Sigma_0,\rho_0)}(y)
\,P(dy)\\
&=
\det(\Gamma_0)
\int
\left(
\mathds{1}_{E(\widetilde{h},(I+\widetilde{A})(I+\widetilde{A})',\rho_0+s)}(z)
-
\mathds{1}_{E(0,I,\rho_0)}(z)
\right)
f(\Gamma_0 z+\mu_0)\,\mathrm{d}z
\\
&=
L'(0,0,0)(\widetilde{h},\widetilde{A},s)
+
o(\|(h,A,s)\|),
\end{split}
\]
by taking $\phi(z)=\det(\Gamma_0)f(\Gamma_0 z+\mu_0)$ in Lemma~\ref{lem:expansion boundary}.
We conclude that
\[
\Lambda_3'(\theta_0)
=
\det(\Gamma_0)
\int_{\partial B_0}
\left(
\frac{\omega'\Gamma_0^{-1}h}{\rho_0} +\frac{\omega'(\Gamma_0^{-1}A+A'\Gamma_0^{-1})\omega}{2\rho_0}+s
\right)
f(\Gamma_0\omega+\mu_0)\sigma_0(d\omega).
\]
For the location functional, with $\theta=(m,G,r)$, we have
\begin{equation}
\label{eq:decomp DLambda1}
\begin{split}
\Lambda_1'(\theta_0)
=
\frac{\partial}{\partial\theta}
\left(
\int_{E(\mu_0,\Sigma_0,\rho_0)}G^{-1}(y-m)\,P(dy)
\right)
\Bigg|_{\theta=\theta_0}
+
\frac{\partial}{\partial\theta}
\left(
\int_{E(m,G^2,r)}\Gamma_0^{-1}(y-\mu_0)\,P(dy)
\right)
\Bigg|_{\theta=\theta_0}.
\end{split}
\end{equation}
The first term on the right hand side of \eqref{eq:decomp DLambda1} can be decomposed as
\begin{equation}
\label{eq:decomp term1}
\frac{\partial}{\partial\theta}
\left(
\int_{E_0}G^{-1}(y-\mu_0)\,P(dy)
\right)
\Bigg|_{\theta=\theta_0}
+
\frac{\partial}{\partial\theta}
\left(
\int_{E_0}\Gamma_0^{-1}(y-m)\,P(dy)
\right)
\Bigg|_{\theta=\theta_0},
\end{equation}
where $E_0=E(\mu_0,\Sigma_0,\rho_0)$.
Because of \eqref{eq:def functionals} and \eqref{eq:phi=ellips},
it follows that
\begin{equation}
\label{eq:property1}
P(E(\mu_0,\Sigma_0,\rho_0))=\gamma
\quad\text{and}\quad
\int_{E(\mu_0,\Sigma_0,\rho_0)}(y-\mu_0)\,P(dy)=0,
\end{equation}
so that the first derivative in \eqref{eq:decomp term1} is equal to zero.
To determine the second derivative in~\eqref{eq:decomp term1},
write
\[
\int_{E(\mu_0,\Sigma_0,\rho_0)}\Gamma_0^{-1}(y-\mu_0-h)\,P(dy)
-
\int_{E(\mu_0,\Sigma_0,\rho_0)}\Gamma_0^{-1}(y-\mu_0)\,P(dy)
=
-\gamma \Gamma_0^{-1}h,
\]
which yields
\[
\frac{\partial}{\partial\theta}
\left(
\int_{E(\mu_0,\Sigma_0,\rho_0)}\Gamma_0^{-1}(y-m)\,P(dy)
\right)
\Bigg|_{\theta=\theta_0}
=
-\gamma \Gamma_0^{-1}h.
\]
For the second term on the right hand of \eqref{eq:decomp DLambda1},
for $(h,A,s)\to(0,0,0)$, consider
\[
\begin{split}
&
\int
\mathds{1}_{E(\mu_0+h,(\Gamma_0+A)(\Gamma_0+A)',\rho_0+s)}(y)\Gamma_0^{-1}(y-\mu_0)
\,P(dy)
-
\int
\mathds{1}_{E(\mu_0,\Sigma_0,\rho_0)}(y)\Gamma_0^{-1}(y-\mu_0)
\,P(dy)\\
&\qquad=
\det(\Gamma_0)
\int
\left(
\mathds{1}_{E(\widetilde{h},(I+\widetilde{A})(I+\widetilde{A})',\rho_0+s)}(z)
-
\mathds{1}_{E(0,I,\rho_0)}(z)
\right)
zf(\Gamma_0 z+\mu_0)\,\mathrm{d}z
\\
&\qquad=
L'(0,0,0)(\widetilde{h},\widetilde{A},s)
+
o(\|(h,A,s)\|),
\end{split}
\]
by taking $\phi(z)=\det(\Gamma_0)zf(\Gamma_0 z+\mu_0)$ in Lemma~\ref{lem:expansion boundary}.
It follows that
\[
\begin{split}
&
\frac{\partial}{\partial\theta}
\left(
\int_{E(m,G^2,r)}\Gamma_0^{-1}(y-\mu_0)\,P(dy)
\right)
\Bigg|_{\theta=\theta_0}\\
&\qquad=
\det(\Gamma_0)
\int_{\partial B_0}
\left(
\frac{\omega'\widetilde{h}}{\rho_0} +\frac{\omega'(\widetilde{A}+\widetilde{A}')\omega}{2\rho_0}+s
\right)
\omega f(\Gamma_0\omega+\mu_0)\sigma_0(d\omega)\\
&\qquad=
\det(\Gamma_0)
\int_{\partial B_0}
\left(
\frac{\omega'\Gamma_0^{-1}h}{\rho_0} +\frac{\omega'(\Gamma_0^{-1}A+A'\Gamma_0^{-1})\omega}{2\rho_0}+s
\right)
\omega f(\Gamma_0\omega+\mu_0)\sigma_0(d\omega).
\end{split}
\]
We conclude that
\[
\Lambda_1'(\theta_0)
=
-\gamma \Gamma_0^{-1}h
+
\det(\Gamma_0)
\int_{\partial B_0}
\left(
\frac{\omega'\Gamma_0^{-1}h}{\rho_0} +\frac{\omega'(\Gamma_0^{-1}A+A'\Gamma_0^{-1})\omega}{2\rho_0}+s
\right)
\omega f(\Gamma_0\omega+\mu_0)\sigma_0(d\omega).
\]
Then, similar to \eqref{eq:decomp DLambda1}, we have
\begin{equation}
\label{eq:decomp DLambda2}
\begin{split}
\Lambda_2'(\theta_0)
&=
\frac{\partial}{\partial\theta}
\left(
\int_{E(\mu_0,\Sigma_0,\rho_0)}
\Big[
G^{-1}(y-m)(y-m)'G^{-1}-I
\Big]\,P(dy)
\right)
\Bigg|_{\theta=\theta_0}\\
&\qquad+
\frac{\partial}{\partial\theta}
\left(
\int_{E(m,G^2,r)}
\Big[
\Gamma_0^{-1}(y-\mu_0)(y-\mu_0)'\Gamma_0^{-1}-I
\Big]\,P(dy)
\right)
\Bigg|_{\theta=\theta_0}.
\end{split}
\end{equation}
The first term in \eqref{eq:decomp DLambda2} can be decomposed as
\begin{equation}
\label{eq:decomp term1 cov}
\begin{split}
&
\frac{\partial}{\partial\theta}
\left(
\int_{E(\mu_0,\Sigma_0,\rho_0)}
\Big[
G^{-1}(y-\mu_0)(y-\mu_0)'\Gamma_0^{-1}-I
\Big]\,P(dy)
\right)
\Bigg|_{\theta=\theta_0}\\
&\quad+
\frac{\partial}{\partial\theta}
\left(
\int_{E(\mu_0,\Sigma_0,\rho_0)}
\Big[
\Gamma_0^{-1}(y-m)(y-\mu_0)'\Gamma_0^{-1}-I
\Big]\,P(dy)
\right)
\Bigg|_{\theta=\theta_0}\\
&\quad+
\frac{\partial}{\partial\theta}
\left(
\int_{E(\mu_0,\Sigma_0,\rho_0)}
\Big[
\Gamma_0^{-1}(y-\mu_0)(y-m)'\Gamma_0^{-1}-I
\Big]\,P(dy)
\right)
\Bigg|_{\theta=\theta_0}\\
&\quad+
\frac{\partial}{\partial\theta}
\left(
\int_{E(\mu_0,\Sigma_0,\rho_0)}
\Big[
\Gamma_0^{-1}(y-\mu_0)(y-\mu_0)'G^{-1}-I
\Big]\,P(dy)
\right)
\Bigg|_{\theta=\theta_0}.
\end{split}
\end{equation}
For the first term in \eqref{eq:decomp term1 cov}, for $\|A\|\to0$, consider
\[
\begin{split}
&
\int_{E(\mu_0,\Sigma_0,\rho_0)}
\Big[
(\Gamma_0+A)^{-1}(y-\mu_0)(y-\mu_0)'\Gamma_0^{-1}-I
\Big]\,P(dy)\\
&\qquad\qquad-
\int_{E(\mu_0,\Sigma_0,\rho_0)}
\Big[
\Gamma_0^{-1}(y-\mu_0)(y-\mu_0)'\Gamma_0^{-1}-I
\Big]\,P(dy)\\
&\quad=
\int_{E(\mu_0,\Sigma_0,\rho_0)}
\Big[
(I-\Gamma_0^{-1}A)\Gamma_0^{-1}(y-\mu_0)(y-\mu_0)'\Gamma_0^{-1}-I
\Big]\,P(dy)
+
o(\|A\|)\\
&\quad=
-\Gamma_0^{-1}A
\int_{E(\mu_0,\Sigma_0,\rho_0)}
\Gamma_0^{-1}(y-\mu_0)(y-\mu_0)'\Gamma_0^{-1}
\,P(dy)+o(\|A\|)
=
-\gamma \Gamma_0^{-1}A+o(\|A\|),
\end{split}
\]
where in the last two steps we use
\begin{equation}
\label{eq:property2}
\int_{E(\mu_0,\Sigma_0,\rho_0)}
\Gamma_0^{-1}(y-\mu_0)(y-\mu_0)'\Gamma_0^{-1}
\,P(dy)=
\gamma I,
\end{equation}
which follows from \eqref{eq:def functionals}. For the second term in \eqref{eq:decomp term1 cov}, consider
\[
\begin{split}
&
\int_{E(\mu_0,\Sigma_0,\rho_0)}
\Big[
\Gamma_0^{-1}(y-\mu_0-h)(y-\mu_0)'\Gamma_0^{-1}-I
\Big]\,P(dy)\\
&\qquad\qquad-
\int_{E(\mu_0,\Sigma_0,\rho_0)}
\Big[
\Gamma_0^{-1}(y-\mu_0)(y-\mu_0)'\Gamma_0^{-1}-I
\Big]\,P(dy)\\
&\quad=
-\Gamma_0^{-1}h
\int_{E(\mu_0,\Sigma_0,\rho_0)}
(y-\mu_0)'\Gamma_0^{-1}
\,P(dy)
=0,
\end{split}
\]
where we use \eqref{eq:property1} and \eqref{eq:property2}.
Because $G$ and $\Gamma_0$ are symmetric, the last two terms in \eqref{eq:decomp term1 cov} are the transpose of the first two terms
in~\eqref{eq:decomp term1 cov}.
This leads to the following derivative for the first term in~\eqref{eq:decomp DLambda2}:
\[
\frac{\partial}{\partial\theta}
\left(
\int_{E(\mu_0,\Sigma_0,\rho_0)}
\Big[
G^{-1}(y-m)(y-m)'G^{-1}-I
\Big]\,P(dy)
\right)
\Bigg|_{\theta=\theta_0}
=
-\gamma (\Gamma_0^{-1}A+A'\Gamma_0^{-1}).
\]
For the second term on the right hand of \eqref{eq:decomp DLambda2},
for $(h,A,s)\to (0,0,0)$, consider
\[
\begin{split}
&
\int
\mathds{1}_{E(\mu_0+h,(\Gamma_0+A)(\Gamma_0+A)',\rho_0+s)}(y)
\left[
\Gamma_0^{-1}(y-\mu_0)(y-\mu_0)'\Gamma_0^{-1}-I
\right]
\,P(dy)\\
&\qquad-
\int
\mathds{1}_{E(\mu_0,\Sigma_0,\rho_0)}(y)
\left[
\Gamma_0^{-1}(y-\mu_0)(y-\mu_0)'\Gamma_0^{-1}-I
\right]
\,P(dy)\\
&\quad=
\det(\Gamma_0)
\int
\left(
\mathds{1}_{E(\widetilde{h},(I+\widetilde{A})(I+\widetilde{A})',\rho_0+s)}(z)
-
\mathds{1}_{E(0,I,\rho_0)}(z)
\right)
\left[
zz'-I
\right]
f(\Gamma_0 z+\mu_0)\,\mathrm{d}z
\\
&\quad=
L'(0,0,0)(\widetilde{h},\widetilde{A},s)
+
o(\|(h,A,s)\|),
\end{split}
\]
by taking $\phi(z)=\det(\Gamma_0)[zz'-I]f(\Gamma_0 z+\mu_0)$ in Lemma~\ref{lem:expansion boundary}.
It follows that
\[
\begin{split}
&
\frac{\partial}{\partial\theta}
\left(
\int_{E(m,G^2,r)}
\Big[
\Gamma_0^{-1}(y-\mu_0)(y-\mu_0)'\Gamma_0^{-1}-I
\Big]\,P(dy)
\right)
\Bigg|_{\theta=\theta_0}
\\
&\qquad=
\det(\Gamma_0)
\int_{\partial B_0}
\left(
\frac{\omega'\widetilde{h}}{\rho_0} +\frac{\omega'(\widetilde{A}+\widetilde{A}')\omega}{2\rho_0}+s
\right)
\left(
\omega\omega'-I
\right)
f(\Gamma_0\omega+\mu_0)\sigma_0(d\omega)\\
&\qquad=
\det(\Gamma_0)
\int_{\partial B_0}
\left(
\frac{\omega'\Gamma_0^{-1}h}{\rho_0} +\frac{\omega'(\Gamma_0^{-1}A+A'\Gamma_0^{-1})\omega}{2\rho_0}+s
\right)
\left(
\omega\omega'-I
\right)
 f(\Gamma_0\omega+\mu_0)\sigma_0(d\omega).
\end{split}
\]
We conclude that
\[
\Lambda_2'(\theta_0)(h,A,s)
=
-\gamma(\Gamma_0^{-1}A+A'\Gamma_0^{-1})
+
\int_{\partial B_0}
\left(
\frac{\omega'\Gamma_0^{-1}h}{\rho_0} +\frac{\omega'(\Gamma_0^{-1}A+A'\Gamma_0^{-1})\omega}{2\rho_0}+s
\right)
\Big(\omega\omega'-I\Big)\, \nu(d\omega).
\]
Noting that we take $A$ symmetric, this finishes the proof.
\hfill\tqed

\bigskip

\noindent
\textbf{Proof of Lemma~\ref{lem:trace=0}:}
From Theorem~\ref{th:derivative},
\begin{eqnarray*}
0 & = & \Lambda_2'(\theta_0)(h,A,s) \\
& = &
-\gamma(\Gamma_0^{-1}A+A\Gamma_0^{-1})
+
\int_{\partial B_0}
\left(
\frac{\omega'\Gamma_0^{-1}h}{\rho_0} +\frac{\omega'(\Gamma_0^{-1}A+A\Gamma_0^{-1})\omega}{2\rho_0}+s
\right)
\Big(\omega\omega'-I\Big) \nu(d\omega),
\end{eqnarray*}
where $B_0=B(0,\rho_0)$.
Taking traces yields
\[
0
=
-2\gamma \mathrm{Tr}(\Gamma_0^{-1}A)
+
(\rho_0^2-k)
\int_{\partial B_0}
\left(
\frac{\omega'\Gamma_0^{-1}h}{\rho_0} +\frac{\omega'(\Gamma_0^{-1}A+A\Gamma_0^{-1})\omega}{2\rho_0}+s
\right)
\,\nu(d\omega).
\]
Because  $\Lambda_3'(\theta_0)=0$, if follows from Theorem~\ref{th:derivative} that the second term on the right hand side
is zero, which proves the lemma.
\hfill\tqed

\bigskip

\noindent
\textbf{Proof of Lemma~\ref{lem:point symmetry}:}
If $f$ satisfies \eqref{eq:point symmetry}, then
\begin{equation}
\label{eq:consequence point symm}
\int_{\partial B_0}\omega_i \omega_j\omega_{m}\,\nu(d\omega)
=
\int_{\partial B_0}\omega_i \,\nu(d\omega)
=
0,
\qquad
\text{for all $i,j,m=1,2,\ldots,k$,}
\end{equation}
where $\nu(d\omega)$ is defined by \eqref{eq:def nu}.
Hence, from Theorem~\ref{th:derivative} we get
\[
0
=
\Lambda_3'(\theta_0)
=
\frac{1}{2\rho_0}
\int_{\partial B_0}
\omega'(\Gamma_0^{-1}A+A\Gamma_0^{-1})\omega\,\nu(d\omega)
+
s\nu(\partial B_0),
\]
which yields the first statement.
Moreover, \eqref{eq:consequence point symm} and Theorem~\ref{th:derivative} also yield that
\[
\begin{split}
0=
\Lambda_1'(\theta_0)
=
-\gamma \Gamma_0^{-1}h
+
\int_{\partial B_0}
\frac{\omega'\Gamma_0^{-1}h}{\rho_0}
\omega \nu(d\omega)
=
\frac1{\rho_0}\left(\int_{\partial B_0}\omega \omega'\,\nu(d\omega)-\gamma \rho_0 I\right)\Gamma_{0}^{-1}h.
\end{split}
\]
Because $f$ satisfies \eqref{eq:point symmetry}, the matrix on the right hand side is a diagonal matrix with elements
\[
\int_{\partial B_0}\omega_i^2\,\nu(d\omega)-\gamma\rho_0.
\]
Therefore, from \eqref{eq:cond h=0} it follows that $h=0$.
\hfill\tqed

\bigskip

The proof of Theorem~\ref{th:non-singular} requires the following lemma.
\begin{lemma}
\label{lem:A=0}
Let $A$ be a $k\times k$ matrix and suppose that $\Gamma A+A\Gamma=0$,
for some $k\times k$ positive definite symmetric matrix $\Gamma$.
Then $A=0$.
\end{lemma}
\textbf{Proof:}
Since $\Gamma$ is positive definite symmetric, there exists a basis of eigenvectors of $\Gamma$. Choose $v$ an eigenvector with eigenvalue $\lambda>0$. Then $\Gamma (Av) + \lambda (Av) = 0$. This means that either $Av$ is an eigenvector of $\Gamma$ with eigenvalue $-\lambda < 0$, which is impossible since $\Gamma$ is positive definite, or $Av=0$. This holds for all eigenvectors of $\Gamma$, and therefore $A=0$.
\hfill\tqed
\bigskip

\noindent
\textbf{Proof of Theorem~\ref{th:non-singular}:}
Suppose that $\Lambda_j'(\theta_0)(h,A,s)=0$ for $j=1,2,3$,
then it suffices to show that $(h,A,s)=(0,0,0)$.
Because \eqref{eq:half space symmetry} implies \eqref{eq:point symmetry}, it follows from
Lemma~\ref{lem:point symmetry} that $h=0$.
Furthermore, $\Lambda_2'(\theta_0)=0$ and $\Lambda_3'(\theta_0)=0$ imply that
\begin{equation}
\label{eq:D2=0}
0=\Lambda_2'(\theta_0)=
-\gamma S
+
\int_{\partial B_0}
\left(
\frac{\omega'S\omega}{2\rho_0}+s
\right)
\omega\omega'\,\nu(d\omega),
\end{equation}
where $B_0=B(0,\rho_0)$ and $S=\Gamma_0^{-1}A+A\Gamma_0^{-1}$ is symmetric.
Condition \eqref{eq:half space symmetry} implies that
\begin{equation}
\label{eq:consequence half space symm}
\int_{\partial B_0}\omega_i \omega_j\omega_{m}\omega_n\,\nu(d\omega)
=
\left\{
\begin{split}
\int_{\partial B_0}\omega_i^2\omega_{m}^2\,\nu(d\omega)&,\quad\text{for }m=n;\,i=j,\\
\int_{\partial B_0}\omega_m^2\omega_{n}^2\,\nu(d\omega)&,\quad\text{for }m\ne n;\,\{i,j\}=\{m,n\},\\
 0&,\quad \text{otherwise,}
\end{split}
\right.
\end{equation}
where $\nu(d\omega)$ is defined by \eqref{eq:def nu}.
Consider the $(m,n)$-th element of equation \eqref{eq:D2=0} for $m\ne n$.
Then it follows from \eqref{eq:consequence half space symm} that
\[
0=
-2\gamma \rho_0 S_{mn}
+
2S_{mn}
\int_{\partial B_0}\omega_m^2\omega_{n}^2\,\nu(d\omega)
=
2\left(
\int_{\partial B_0}\omega_m^2\omega_{n}^2\,\nu(d\omega)-\gamma \rho_0
\right)
S_{mn}.
\]
The factor in front of $S_{mn}$ is non-zero by assumption~\eqref{eq:cond Anm=0}, so that $S_{mn}=0$
for all $m\neq n$.
Finally, consider the $(m,m)$-th element of~\eqref{eq:D2=0}
and insert \eqref{eq:relation s-A}, which is obtained from Lemma~\ref{lem:point symmetry}.
Then we get
\[
0
=
-2\gamma \rho_0 S_{mm}
+
\sum_{i=1}^k
\left(
\int_{\partial B_0}\omega_i^2\omega_{m}^2\,\nu(d\omega)
-
\frac1{\nu_0}
\int_{\partial B_0}\omega_i^2\,\nu(d\omega)
\int_{\partial B_0}\omega_m^2\,\nu(d\omega)
\right)
S_{ii}.
\]
The right hand side is of the form $Mx$,
where $x=\mathrm{diag}(S)$ and $M$ is defined in \eqref{eq:def M}.
However, since $\mathrm{Tr}(S)=0$ according to Lemma~\ref{lem:trace=0},
from \eqref{eq:property M} we conclude $S_{mm}=0$ for all $m=1,2,\ldots,k$.
It follows that $S=\Gamma_0^{-1}A+A\Gamma_0^{-1}=0$, and consequently, by~\eqref{eq:relation s-A}, we have $s=0$.
Furthermore, from Lemma~\ref{lem:A=0} we conclude that $A=0$.
\hfill\tqed
\subsection{Proofs for Section~\ref{sec:elliptical contoured}}
\textbf{Proof of Proposition~\ref{prop:cond elliptical}:}
Conditions (B) and \eqref{eq:half space symmetry} are immediate.
From \eqref{eq:MCD at ellip cont}, we find that
$f(\Gamma_0\omega+\mu_0)$ is constant on $\partial B_0$:
\[
f(\Gamma_0\omega+\mu_0)
=
\det(\Sigma)^{-1/2}
h(\alpha(\gamma)^2\|\omega\|^2)
=
\alpha(\gamma)^k\det(\Gamma_0)^{-1}
h(r(\gamma)^2),
\]
and
\begin{equation}
\label{eq:def nu ellip cont}
\nu(d\omega)
=
\alpha(\gamma)^kh(r(\gamma)^2)\sigma_0(d\omega),
\end{equation}
for $\omega\in \partial B_0$.
One can easily check that for all $i,j=1,2,\ldots,k$ (e.g., see Lemma 1 in~\cite{lopuhaa97})
\begin{equation}
\label{eq:integrals}
\begin{split}
\int_{\partial B_0}\omega_i^2\,\nu(d\omega)
&=
\frac1k\int_{\partial B_0}\|\omega\|^2\,\nu(d\omega)
=
\frac{2\pi^{k/2}}{k\Gamma(k/2)}h(r(\gamma)^2)r(\gamma)^k\rho_0,\\
\int_{\partial B_0}\omega_i^2\omega_{j}^2\,\nu(d\omega)
&=
\frac{1+2\delta_{ij}}{k(k+2)}
\int_{\partial B_0}\|\omega\|^4\,\nu(d\omega)\\
&=
(1+2\delta_{ij})
\frac{2\pi^{k/2}}{k(k+2)\Gamma(k/2)}h(r(\gamma)^2) r(\gamma)^k\rho_0^{3},\\
\nu(\partial B_0)
&=
\frac{2\pi^{k/2}}{\Gamma(k/2)}h(r(\gamma)^2) r(\gamma)^{k-1}\alpha(\gamma)
=
\frac{k}{\rho_0^2}\int_{\partial B_0}\omega_i^2\,\nu(d\omega).
\end{split}
\end{equation}
Because $h$ is decreasing
and non-constant on $[0,r(\gamma))$, conditions \eqref{eq:cond h=0} and \eqref{eq:cond Anm=0} are fulfilled:
\begin{equation}
\label{eq:cond nonsingular ellip cont}
\begin{split}
\gamma
&=
\frac{2\pi^{k/2}}{\Gamma(k/2)}
\int_0^{r(\gamma)}
h(r^2) r^{k-1}\,dr
>
\frac{2\pi^{k/2}}{k\Gamma(k/2)}h(r(\gamma)^2) r(\gamma)^{k},\\
\gamma\alpha(\gamma)^2
&=
\frac{2\pi^{k/2}}{k\Gamma(k/2)}\int_0^{r(\gamma)}h(r^2)r^{k+1}\,dr
>
\frac{2\pi^{k/2}}{k(k+2)\Gamma(k/2)}h(r(\gamma)^2)r(\gamma)^{k+2}.
\end{split}
\end{equation}
Finally, from the equations above, it follows that the matrix $M$ defined in~\eqref{eq:def M} can be decomposed as
$M=c_1I+c_2\mathbf{11}'$, where $\mathbf{1}=(1,1,\ldots,1)'$,
and
\[
\begin{split}
c_1
&=
\frac{4\pi^{k/2}}{k(k+2)\Gamma(k/2)}h(r(\gamma)^2)r(\gamma)^k \rho_0^{3}-2\gamma \rho_0,\\
c_2
&=
-\frac{4\pi^{k/2}}{k^2(k+2)\Gamma(k/2)}h(r(\gamma)^2) r(\gamma)^k\rho_0^{3}.
\end{split}
\]
Because
$Mx=c_1x+c_2(x_1+\cdots+x_k)\mathbf{1}$,
it follows that, if $x_1+\cdots+x_k=0$,
$Mx=0$ implies $x=0$ as long as $c_1\ne 0$,
i.e.,
\[
\frac{2\pi^{k/2}}{k(k+2)\Gamma(k/2)}h(r(\gamma)^2) r(\gamma)^{k+2}\ne \gamma\alpha(\gamma)^{2},
\]
which follows from \eqref{eq:cond nonsingular ellip cont}.
\hfill\tqed

\bigskip

\noindent
\textbf{Proof of Theorem~\ref{th:derivative ellip}:}
From~\eqref{eq:def nu ellip cont}, it follows that
\begin{equation}
\label{eq:consequence1}
\begin{split}
\int_{\partial B_0}\omega_i \omega_j\,\nu(d\omega)&=0,
\quad
\text{ for }i\ne j\\
\int_{\partial B_0}\omega_i\,\nu(d\omega)=0
\quad\text{and}\quad
\int_{\partial B_0}\omega_i \omega_j\omega_k\,\nu(d\omega)&=0,
\quad
\text{ for all }i, j, k,
\end{split}
\end{equation}
Hence, from Theorem~\ref{th:derivative} we find
\[
\Lambda_1'(\theta_0):(h,A,s)
=
-\gamma \Gamma_0^{-1}h
+
\frac1{\rho_0}
\int_{\partial B_0}
\omega\omega'\Gamma_0^{-1}h
\, \nu(d\omega)
=
\beta_1h,
\]
where according to \eqref{eq:integrals} and \eqref{eq:cond nonsingular ellip cont},
\[
\beta_1
=
\frac1{k\alpha \rho_0}
\int_{\partial B_0}
\|\omega\|^2
\, \nu(d\omega)
-
\frac{\gamma}{\alpha }
=
\dfrac1{\alpha}
\left(
\dfrac{\rho_0}{k}\nu_0-\gamma
\right)<0.
\]
Next, consider $\Lambda_3'(\theta_0)$.
From~\eqref{eq:consequence1} we find
\[
\Lambda_3'(\theta_0)
=
\frac1{\rho_0}\int_{\partial B_0}\omega'\Gamma_0^{-1}A\omega\,\nu(d\omega)+s\nu_0.
\]
From~\eqref{eq:integrals}, the first term on the right hand side is
\[
\frac1{\alpha \rho_0}\int_{\partial B_0}\omega'A\omega\,\nu(d\omega)
=
\frac1{k\alpha \rho_0}\int_{\partial B_0}\|\omega\|^2\mathrm{Tr}(A)\,\nu(d\omega)
=
\frac{\rho_0\nu_0}{k\alpha }\mathrm{Tr}(A).
\]
This means that $\Lambda_3'(\theta_0)=\beta_5\mathrm{Tr}(A)+\beta_6s$.
Finally, from~\eqref{eq:MCD at ellip cont} and~\eqref{eq:consequence1},
\[
\begin{split}
\Lambda_2'(\theta_0)
&=
-\gamma(\Gamma_0^{-1}A+A\Gamma_0^{-1}) +
\int_{\partial B_0}
\left(
\frac{\omega'(\Gamma_0^{-1}A+A\Gamma_0^{-1})\omega}{2\rho_0}+s
\right)
\Big(\omega\omega'-I\Big)\, \nu(d\omega)\\
&=
-\frac{2\gamma}{\alpha }A
-
\Lambda_3'(\theta_0)\cdot I
+
\frac{1}{\alpha \rho_0}
\int_{\partial B_0}
(\omega'A\omega)\omega\omega'\, \nu(d\omega)
+
s\int_{\partial B_0}
\omega\omega'\, \nu(d\omega)\\
&=
-\frac{2\gamma}{\alpha }A
-
(\beta_5\mathrm{Tr}(A)+\beta_6s)\cdot I
+
\frac{1}{\alpha \rho_0}
\int_{\partial B_0}
(\omega'A\omega)\omega\omega'\, \nu(d\omega)
+
\frac{\rho_0^2}{k}\nu_0\cdot sI.
\end{split}
\]
Consider the $(m,n)$-th element of the third integral on the right hand side.
From~\eqref{eq:consequence half space symm} and~\eqref{eq:integrals}, it follows that this integral is equal to
\[
\frac{1}{\alpha \rho_0}
\sum_{i=1}^k\sum_{j=1}^k A_{ij}
\int_{\partial B_0}
\omega_i\omega_j\omega_m\omega_n\, \nu(d\omega)
=
\frac{\rho_0^3\nu_0}{\alpha k(k+2)}
\left(
\mathrm{Tr}(A)\mathds{1}_{\{m=n\}}+2A_{mn}
\right),
\]
which means that
\[
\frac{1}{\alpha \rho_0}
\int_{\partial B_0}
(\omega'A\omega)\omega\omega'\, \nu(d\omega)
=
\frac{\rho_0^3\nu_0}{\alpha k(k+2)}
\left(
\mathrm{Tr}(A)\cdot I+2A
\right).
\]
Summarizing, in the expression of $\Lambda_2'(\theta_0)$, the coefficient of $A$ is
\[
\beta_2
=
\frac{2\rho_0^3\nu_0}{\alpha k(k+2)}
-
\frac{2\gamma}{\alpha },
\]
the coefficient of $\mathrm{Tr}(A)\cdot I$ is
\[
\beta_3
=
\frac{\rho_0^3\nu_0}{\alpha k(k+2)}
-
\frac{\rho_0\nu_0}{k\alpha },
\]
and the coefficient of $sI$ is
\[
\beta_4
=\frac{\rho_0^2}{k}\nu_0
-
\nu_0.
\]
From \eqref{eq:integrals} and \eqref{eq:cond nonsingular ellip cont}, it can be seen that $\beta_2<0$.

To determine the expression of the inverse mapping, put $D(h,A,s)=(g,B,t)$ and solve for $(h,A,s)$.
For the vector valued component of $D$, we have $g=D_1(h,A,s)=\beta_1h$.
Since $\beta_1<0$, this immediately gives $h=\beta_1^{-1}g$.
For the remaining mappings put
\begin{equation}
\label{eq:D2D3}
\begin{split}
B&=D_2(h,A,s)=\beta_2A+\beta_3\mathrm{Tr}(A)\cdot I+\beta_4s\cdot I\\
t&=D_3(h,A,s)=\beta_5\mathrm{Tr}(A)+\beta_6s.
\end{split}
\end{equation}
By taking traces in the first equation we can solve for $\mathrm{Tr}(A)$ and $s$:
\begin{equation}
\label{eq:TrA and s}
\begin{split}
c\mathrm{Tr}(A)&=\beta_6 \mathrm{Tr}(B)-k\beta_4t\\
cs&=(\beta_2+k\beta_3)t-\beta_5\mathrm{Tr}(B),
\end{split}
\end{equation}
where
$c=\beta_2\beta_6+k\beta_3\beta_6-k\beta_4\beta_5=-2\gamma\beta_6/\alpha$.
Since $\beta_2<0$ and $\beta_6>0$, from \eqref{eq:D2D3} and \eqref{eq:TrA and s} it follows that
\[
\begin{split}
A
&=
\beta_2^{-1}
\left(
B-\beta_3\mathrm{Tr}(A)\cdot I-\beta_4s\cdot I
\right)\\
&=
\beta_2^{-1}B
-\frac{\beta_3}{c\beta_2}\left(\beta_6\mathrm{Tr}(B)-k\beta_4t\right)\cdot I
-\frac{\beta_4}{c\beta_2}\left(-\beta_5\mathrm{Tr}(B)+(\beta_2+k\beta_3)t\right)\cdot I\\
&=
\beta_2^{-1}B
+\frac{\alpha(\beta_3\beta_6-\beta_4\beta_5)}{2\gamma\beta_2\beta_6}\mathrm{Tr}(B)\cdot I
-\frac{\alpha\beta_2\beta_3}{2\gamma\beta_2\beta_6}t\cdot I
\end{split}
\]
and
\[
\begin{split}
s&=-\frac{\beta_5}c\mathrm{Tr}(B)+\frac{\beta_2+k\beta_3}ct
=
\frac{\alpha\beta_5}{2\gamma\beta_6}\mathrm{Tr}(B)-\frac{\alpha(\beta_2+k\beta_3)}{2\gamma\beta_6}t.
\qquad
\hfill\tqed
\end{split}\]

\bigskip

\noindent
\textbf{Proof of Corollary~\ref{cor:expansion}:}
Since $\Lambda(\theta_0)^{-1}$ is a linear mapping
and $\mathds{E}\Psi(X_i,\theta_0)=0$, we obtain from~\eqref{eq:expansion estimators},
\begin{equation}
\label{eq:expansion2}
\widehat{\theta}_n-\theta_0
=
-\frac1n\sum_{i=1}^n
\Lambda'(\theta_0)^{-1}\Psi(X_i,\theta_0)+o_\mathbb{P}(n^{-1/2}),
\end{equation}
where $\Psi$ is defined in~\eqref{eq:def Psi}.
In particular, we have
\begin{equation}
\label{eq:expansion mu,rho}
\begin{split}
\sqrt{n}\widehat{\mu}_n
&=
-\frac1{\sqrt{n}}\sum_{i=1}^n
\left[D^{\mathrm{inv}}\right]_1\Psi(X_i,\theta_0)+o_\mathbb{P}(1),\\
\sqrt{n}
\left(
\widehat{\rho}_n-\frac{r}{\alpha}
\right)
&=
-\frac1{\sqrt{n}}\sum_{i=1}^n
\left[D^{\mathrm{inv}}\right]_3\Psi(X_i,\theta_0)+o_\mathbb{P}(1).
\end{split}
\end{equation}
According to \eqref{eq:MCD at ellip cont}, $\theta_0=(\mu_0,\Gamma_0,\rho_0)=(0,\alpha I,r/\alpha)$, so that
\begin{equation}
\label{eq:Psi0}
\begin{split}
\Psi_1(x,\theta_0)&=\mathds{1}_{\{\|x\|\leq r\}}\alpha^{-1}x,\\
\Psi_2(x,\theta_0)&=\mathds{1}_{\{\|x\|\leq r\}}\left(\alpha^{-2}xx'-I\right),\\
\Psi_3(x,\theta_0)&=\mathds{1}_{\{\|x\|\leq r\}}-\gamma.
\end{split}
\end{equation}
Insert $g_0=\Psi_1(x,\theta_0)$, $B_0=\Psi_2(x,\theta_0)$, and $t_0=\Psi_3(x,\theta_0)$
in the expressions for $D^{\mathrm{inv}}(g,B,t)$ given in Theorem~\ref{th:derivative ellip}.
Then we find
\begin{equation}
\label{eq:expression Dinv13}
\begin{split}
\left[D^{\mathrm{inv}}\right]_1\Psi(x,\theta_0)
&=
(\alpha\beta_1)^{-1}\mathds{1}_{\{\|x\|\leq r\}}x\\
\left[D^{\mathrm{inv}}\right]_3\Psi(x,\theta_0)
&=
\frac{\alpha\beta_5}{2\gamma\beta_6}
\mathds{1}_{\{\|x\|\leq r\}}\left(\frac{\|x\|^2}{\alpha^2}-k\right)
-\frac{\alpha(\beta_2+k\beta_3)}{2\gamma\beta_6}
\left(
\mathds{1}_{\{\|x\|\leq r\}}-\gamma
\right).
\end{split}
\end{equation}
Together with \eqref{eq:expansion mu,rho}, this immediately yields the expansion for $\sqrt{n}\widehat{\mu}_n$ and
the expansion for $\sqrt{n}\left(\widehat{\rho}_n-r/\alpha\right)$ with
\[
\begin{split}
\lambda_1
&=
-\frac{\beta_5}{2\alpha\gamma\beta_6}
=
-\frac{r}{2k\gamma\alpha^3},\\
\lambda_2
&=
\frac{\alpha(\beta_2+k\beta_3+k\beta_5)}{2\gamma\beta_6}
=
\frac{r^3}{2k\gamma\alpha^3}-\frac{1}{\beta_6},\\
\lambda_3
&=
-\frac{\alpha(\beta_2+k\beta_3)}{2\beta_6}
=
\frac{\gamma}{\beta_6}+\frac{r}{2k\alpha^3}\left(k\alpha^2-r^2\right).
\end{split}
\]
To obtain the expansion for the covariance estimator, note that $P$ satisfies the conditions
of Theorem~4.2 in~\cite{catorlopuhaa2009}.
This means $\widehat{\Gamma}_n\to\alpha I$ with probability one, so that
\[
\widehat{\Sigma}_n-\alpha^2 I
=
(\widehat{\Gamma}_n+\alpha I)(\widehat{\Gamma}_n-\alpha I)
=
2\alpha
(\widehat{\Gamma}_n-\alpha I)
+o(1),
\]
with probability one.
Hence, from \eqref{eq:expansion2} we obtain
\begin{equation}
\label{eq:expansion Sigma}
\sqrt{n}
\left(
\widehat{\Sigma}_n-\alpha^2 I
\right)
=
-\frac{2\alpha}{\sqrt{n}}
\sum_{i=1}^n
\left[D^{\mathrm{inv}}\right]_2\Psi(X_i,\theta_0)+o_\mathds{P}(1),
\end{equation}
where
\begin{equation}
\label{eq:expression Dinv2}
\begin{split}
\left[D^{\mathrm{inv}}\right]_2\Psi(x,\theta_0)
&=
\beta_2^{-1}
\mathds{1}_{\{\|x\|\leq r\}}\left(\frac{xx'}{\alpha^2}-I\right)\\
&\qquad+
\frac{\alpha(\beta_3\beta_6-\beta_4\beta_5)}{2\gamma\beta_2\beta_6}
\mathds{1}_{\{\|x\|\leq r\}}\left(\frac{\|x\|^2}{\alpha^2}-k\right)\cdot I\\
&\qquad+
\frac{\alpha\beta_4}{2\gamma\beta_6}
\left(\mathds{1}_{\{\|x\|\leq r\}}-\gamma\right)\cdot I.
\end{split}
\end{equation}
This yields the expansion for $\sqrt{n}\left(\widehat{\Sigma}_n-\alpha^2 I\right)$ with
\[
\begin{split}
\kappa_1
&=
\frac{2\alpha}{\beta_2}
+\frac{k\alpha^2(\beta_3\beta_6-\beta_4\beta_5)}{\gamma\beta_2\beta_6}
-\frac{\alpha^2\beta_4}{\gamma\beta_6}
=
-\frac{r^2}{k\gamma},\\
\kappa_2
&=
\frac{\beta_4\beta_5-\beta_3\beta_6}{\gamma\beta_2\beta_6}
=
\dfrac{\alpha\beta_2+2\gamma}{k\gamma\alpha\beta_2},
\\
\kappa_3
&=
-\frac{2}{\alpha\beta_2},\\
\kappa_4
&=
\frac{\alpha^2\beta_4}{\beta_6}
=
\frac{r^2-k\alpha^2}{k}.
\qquad\hfill\tqed
\end{split}
\]

\bigskip

\noindent
\textbf{Proof of Theorem~\ref{th:AN}:}
The expansion for $\sqrt{n}\widehat{\mu}_n$ given in Corollary~\ref{cor:expansion},
together with the fact that $\mathbb{E}\{\|X_1\|\leq r\}X_1=0$,
yields that $\sqrt{n}\widehat{\mu}_n$ is asymptotically normal with mean zero and covariance matrix
\[
\pi^2
\mathbb{E}
\{\|X_1\|\leq r\}X_1X_1'
=
\frac{\pi^2}k
\mathbb{E}\{\|X_1\|\leq r\}\|X_1\|^2\cdot I.
\]
Since $\pi=-(\alpha\beta_1)^{-1}$, together with \eqref{eq:def alpha}, we find
\[
\tau
=
\frac{\pi^2}k
\mathbb{E}\{\|X_1\|\leq r\}\|X_1\|^2
=
\frac{k^2\gamma\alpha^4}{(k\gamma\alpha-r\nu_0)^2},
\]
which proves part (ii).
To prove (iii), first note that from Corollary~\ref{cor:expansion}, it follows that
\begin{equation}
\label{eq:expansion Sigma ell m}
\sqrt{n}
(\widehat{\Sigma}_n-\alpha^2 I)
=
\frac{1}{\sqrt{n}}
\sum_{i=1}^n
\left(
\ell(\|X_i\|)\frac{X_iX_i'}{\|X_i\|^2}+m(\|X_i\|)\cdot I
\right)
+o_\mathbb{P}(1)
\end{equation}
where
$\ell(y)=\kappa_3\mathds{1}_{\{\|y\|\leq r\}}y^2$
and $m(y)=\mathds{1}_{\{\|y\|\leq r\}}(\kappa_1+\kappa_2y^2)+\kappa_4$.
Note that according to~\eqref{eq:def alpha},
\[
\begin{split}
\mathbb{E}\ell(\|X_1\|)
&=
-\frac{2}{\alpha\beta_2}\mathbb{E}\{\|X_1\|\leq r\}\|X_1\|^2
=
-\frac{2\alpha k\gamma}{\beta_2},\\
\mathbb{E}m(\|X_1\|)
&=
-\frac{r^2}{k\gamma}
\mathbb{E}\{\|X_1\|\leq r\}
+
\frac{\alpha\beta_2+2\gamma}{k\gamma\alpha\beta_2}
\mathbb{E}\{\|X_1\|\leq r\}\|X_1\|^2
+
\frac{r^2-k\alpha^2}{k}
=
\frac{2\alpha\gamma}{\beta_2},
\end{split}\]
so that $\mathbb{E}\left[\ell(\|X_1\|)+km(\|X_1\|)\right]=0$.
Since also $\mathbb{E} \ell^2(\|X_1\|)<\infty$ and $\mathbb{E} m^2(\|X_1\|)<\infty$,
it follows from Lemma~5 in~\cite{lopuhaa97}, that
the sum on the right hand side of~\eqref{eq:expansion Sigma ell m}
is asymptotically normal with mean zero and covariance matrix
$\sigma_1(I+C_{k,k})+\sigma_2\mbox{vec}(I)\mbox{vec}(I)'$,
where
\[\begin{split}
\sigma_1
&=
\frac{\mathbb{E}\,\ell^2(\|X_1\|)}{k(k+2)},
\\
\sigma_2 &=
\frac{ \mathbb{E}\,\ell^2(\|X_1\|)}{k(k+2)}
            +  \mathbb{E}\,m^2(\|X_1\|)
            + \frac2k \mathbb{E}\,\ell(\|X_1\|)m(\|X_1\|).
\end{split}\]
If we fill in the expressions for $\ell(\|X_1\|$ and $m(\|X_1\|)$, we get
\[
\begin{split}
\sigma_1&=\dfrac{\kappa_3^2}{k(k+2)}\mathbb{E}\mathds{1}_{\{\|X_1\|\leq r\}}\|X_1\|^4\\
\sigma_2
&=
\left(
\frac{\kappa_3^2}{k(k+2)}+\kappa_2^2+\frac{2\kappa_2\kappa_3}{k}
\right)
\mathbb{E}\mathds{1}_{\{\|X_1\|\leq r\}}\|X_1\|^4\\
&\quad+
\left(
2(\kappa_1+\kappa_4)\kappa_2+\frac2k\kappa_3(\kappa_1+\kappa_4)
\right)
\mathbb{E}\mathds{1}_{\{\|X_1\|\leq r\}}\|X_1\|^2
+
\kappa_1(\kappa_1+2\kappa_4)
\mathbb{E}\mathds{1}_{\{\|X_1\|\leq r\}}
+
\kappa_4^2.
\end{split}\]
Substituting the expressions for $\kappa_1,\kappa_2,\kappa_4$ given in Corollary~\ref{cor:IF}
together with \eqref{eq:def alpha} and \eqref{eq:def r(gamma)} proves (iii).
For part (iv) note that
\[
\mathbb{E}
\left[
\lambda_1\mathds{1}_{\{\|X_i\|\leq r\}}\|X_i\|^2
+
\lambda_2\mathds{1}_{\{\|X_i\|\leq r\}}
+
\lambda_3
\right]
=
\lambda_1k\gamma\alpha^2
+
\lambda_2\gamma
+
\lambda_3
=0.
\]
Therefore, from the expansion given in Corollary~\ref{cor:expansion}, it follows that $\sqrt{n}(\widehat{\rho}_n-r/\alpha)$
is asymptotically normal with variance
\[
\begin{split}
\sigma_\rho^2
&=
\mathbb{E}
\left(
\lambda_1\mathds{1}_{\{\|X_i\|\leq r\}}\|X_i\|^2
+
\lambda_2\mathds{1}_{\{\|X_i\|\leq r\}}
+
\lambda_3
\right)^2\\
&=
\lambda_1^2\mathbb{E}\mathds{1}_{\{\|X_i\|\leq r\}}\|X_i\|^4
+
\lambda_1(\lambda_2+\lambda_3)
\mathbb{E}\mathds{1}_{\{\|X_i\|\leq r\}}\|X_i\|^2
+
\lambda_2(\lambda_2+\lambda_3)\mathbb{E}\mathds{1}_{\{\|X_i\|\leq r\}}
+
\lambda_3^2.
\end{split}
\]
Substituting the expressions for $\lambda_2,\lambda_3$ given in Corollary~\ref{cor:IF}
together with \eqref{eq:def alpha} and \eqref{eq:def r(gamma)} proves (iv).
Finally, for part (i), first note that according to Theorem~5.1 in~\cite{catorlopuhaa2009},
$\widehat{\mu}_n$, $\widehat{\Sigma}_n$ and $\widehat{\rho}_n$ are mutually asymptotically normal.
Hence, it suffices to prove that the quantities considered in part (i) are asymptotically uncorrelated.
However, this follows directly from the expansions given in Corollary~\ref{cor:expansion}
together with the symmetry properties of spherically symmetric densities.
\hfill\tqed

\bigskip

\noindent
\textbf{Proof of Corollary~\ref{cor:IF}:}
According to Theorem~1 in \cite{butlerdaviesjuhn93}, the MCD functional $\theta_0=(\mu_0,\Gamma_0,\rho_0)$ as defined in \eqref{eq:def theta} is unique,
and since $P$ has a density, all conditions of Theorem~5.2 in \cite{catorlopuhaa2009} are satisfied.
It follows from this theorem that the influence function
for the functional $\Theta(P)=\big(\mu(P),\Gamma(P),\rho(P)\big)$, where $\Gamma(P)^2=\Sigma(P)$,
is given by
\begin{equation}
\label{eq:general expression IF}
\mathrm{IF}(x;\Theta,P)=-\Lambda'(\theta_0)^{-1}\Psi(x,\theta_0),
\end{equation}
where $\Psi$ is defined in~\eqref{eq:def Psi}.
The expressions for $\mathrm{IF}(x;\mu,P)$ and $\mathrm{IF}(x;\rho,P)$
follow directly from~\eqref{eq:expression Dinv13}.
To obtain the influence function for the covariance functional, first note that
according to the continuity of the MCD functional,
$\Gamma(P_{\eps,x})\to\Gamma(P)=\alpha I$, as $\eps\downarrow0$, where $P_{\eps,x}=(1-\eps)P+\eps\delta_x$.
This means that
\[
\Sigma(P_{\eps,x})-\Sigma(P)
=
\big(
\Gamma(P_{\eps,x})+\Gamma(P)
\big)
\big(
\Gamma(P_{\eps,x}-\Gamma(P)
\big)
=
2\alpha
\big(
\Gamma(P_{\eps,x})-\Gamma(P)
\big)
+o(\eps).
\]
It follows that
\[
\mathrm{IF}(x;\Sigma,P)
=
2\alpha\cdot
\mathrm{IF}(x;\Gamma,P)
=
-2\alpha \left[D^{\mathrm{inv}}\right]_2\Psi(x,\theta_0).
\]
The expression then follows from~\eqref{eq:expression Dinv2}.
\hfill\tqed

\end{document}